\documentclass[11pt]{article}
\usepackage{amsmath,amsthm,amsfonts,latexsym}

\begin{document}

\def\qed{\vbox{\hrule
  \hbox{\vrule\hbox to 5pt{\vbox to 8pt{\vfil}\hfil}\vrule}\hrule}}
 
\def\endproof{\unskip \nobreak \hskip0pt plus 1fill \qquad \qed \par}

\newcommand{\N}{{\mathbb{N}}}
\newcommand{\R}{{\mathbb{R}}}
\newcommand{\C}{{\mathbb{C}}}
\newcommand{\T}{{\mathbb{T}}}
\newcommand{\Z}{{\mathbb{Z}}}
\newcommand{\F}{{\mathbb{F}}}

\def\B{\mathcal{B}}

\def\P{\mathbb{P}}
\def\H{\mathcal{H}}
\def\K{\mathcal{K}}
\newcommand{\HH}{\mathfrak H}
\newcommand{\KK}{\mathfrak K}
\newcommand{\LL}{\mathfrak L}

\newcommand{\no}[1]{|| {#1}Ê||}

\newcommand{\ve}[1]{\mathbf{#1}}
\newcommand{\mve}[1]{$\mathbf{#1}$}

\title{The full group C$^*$-algebra of the modular group is primitive}

\author{Erik B\'edos$^{*}$,
Tron {\AA}. Omland$^{*}$ \\
\\
\\}

\date{{\it Revised Mars 29, 2010}}
\maketitle
\renewcommand{\sectionmark}[1]{}

\begin{abstract}
We show that the full group C$^*$-algebra of $PSL(n, \Z)$ is primitive when $n=2$, and not primitive when $n\geq 3$. Moreover, we show that there exists an uncountable family of pairwise inequivalent, faithful irreducible representations of $C^*(PSL(2,\Z))$.

\medskip  \medskip \medskip  
\vskip 0.9cm

 {\it \qquad Dedicated to the memory of Gerard J.\  Murphy.}

\medskip  \medskip \medskip 
\vskip 0.9cm
\noindent {\bf MSC 2010}:  Primary: 46L05. Secondary: 22D25, 46L55.

\smallskip
\noindent {\bf Keywords}: 
 modular group,  full group $C^{*}$-algebra, primitivity, twisted crossed product
\end{abstract}

\vfill
\thanks{\noindent $^*$ partially supported by the Norwegian Research
Council.\\

\par \par
 
 \noindent E.B.{'}s address: Institute of Mathematics, University of
Oslo, P.B. 1053 Blindern, 0316 Oslo, Norway. E-mail : bedos@math.uio.no. \\

\noindent 
T.O.{'}s address: Department of Mathematical Sciences, NTNU,
7491 Trondheim, Norway. E-mail : tronanen@math.ntnu.no.\par}

\newpage
\begin{flushleft}
\section{Introduction} 

Simple and, more generally, primitive and prime C$^*$-algebras may be considered as building blocks of the theory, playing a somewhat similar role as factors do within the theory of von Neumann algebras.  If we restrict ourselves to separable C$^*$-algebras, as we always do in this paper, primitivity is equivalent to primeness (see for example \cite{Ped}), and we will therefore refer to primitivity for both notions.  Now, given some class of separable C$^*$-algebras, one natural task is to investigate which members of this class are simple or primitive.

\medskip An interesting family of separable C$^*$-algebras consists of the group C$^*$-algebras associated with countable discrete groups.    
We recall that such a group $G$ is called  C${^*}$-simple if its  {\it reduced} group C$^*$-algebra $C_{r}^{*}(G)$ is simple. As the  {\it full} group C$^*$-algebra $C^{*}(G)$ is simple only when $G$ is trivial, this terminology is not ambiguous.  
The class of C$^{*}$-simple groups has received a lot of attention during the last decades and the reader may consult \cite{dHa} for a recent, comprehensive review.
It is also well known (see \cite{Pa, Mur})
that $C_{r}^{*}(G)$ is primitive if and only if $G$ is icc (that is, every nontrivial conjugacy class in 
$G$ is infinite), if and only if the group von Neumann algebra of $G$ is a factor.  f

\medskip On the other hand, the problem of determining when $C^{*}(G)$ is primitive seems hard in general. 
A necessary condition is that $G$ is icc \cite{Mur}, and this
condition is  also sufficient when $G$ is assumed to be amenable, as $C^*(G)$ is then isomorphic to $C^*_{r}(G)$.
We note in passing that this problem is quite different from the one of determining the class of  groups having a faithful irreducible unitary representation, which contains many other groups besides all icc groups (see \cite{BH2}).

\medskip Until a few years ago, the only known nonamenable icc 
 groups having a primitive full group C$^*$-algebra were nonabelian free groups, as originally shown by  H.\ Yoshizawa \cite{Y} 
 and rediscovered later by M.D.\ Choi  \cite{Cho, Da}. Then primitivity of $C^{*}(G)$ was established when $G=G_{1}*G_{2}$ is the  free product of two countable subgroups $G_1$ and $G_2$  satisfying at least one of the following assumptions:  
  \begin{itemize}
  \item[i)] $G_1=\Z*\Z$ or $G_1= \Z *\Z_2$ ($G_2$ being any  group).
 \item[ii)]  $G_1$ is nontrivial and free, and $G_2$ is a nontrivial and amenable.
  \item[iii)] 
  $G_1$ is a nonabelian and free, and $C^{*}(G_2)$ admits no nontrivial projections.
 \end{itemize} 
Case (i) is due to N.\ Khatthou \cite[Th\' eor\`emes 2 et 3]{Kha},  while (ii) and (iii) are due to G.J.\ Murphy 
  \cite[Theorems 3.3 and 3.4]{Mur}.

\medskip 
In \cite[Problem 25]{dHa}, P.\ de la Harpe raises the problem of finding other (nonamenable icc) groups having a primitive 
full group C$^{*}$-algebra. One may especially wonder whether this property holds for any group $G$ which is the free product  of two nontrivial groups, where at least one of them has more than two elements (as the infinite dihedral group $\Z_2 *\Z_2$ is not icc). The simplest case for which the answer is unknown  is  the modular group  $PSL(2,\Z)=$ $ \Z_{2}*\Z_{3}$,
and our main result in this paper is  that  $C^*(PSL(2,\Z))$ is indeed primitive (cf.\ Theorem 2).

\medskip An outline of our proof is as follows. 
Let $H$ be the kernel of the canonical homomorphism from $G=  \Z_{2}*\Z_{3}$ onto
$ \Z_{2}\times \Z_{3}$. Then $H$ $\simeq  \Z*\Z$. 
Exploiting a  certain phase-action of the circle group $\T$ on $C^*(H)$, we then show how a faithful irreducible representation of $C^*(H)$  
may be picked so that it induces a representation of $C^*(G)$ which is also faithful and irreducible. Moreover, we show that there exists an uncountable family of  pairwise inequivalent, irreducible faithful representations of $C^*(G)$.  \\ A similar idea was used by Murphy in his proof of  \cite[Theorem 3.3]{Mur}, where he considers certain semidirect products of nonabelian free groups by amenable groups. However, in our case,  the exact sequence $1 \rightarrow H \rightarrow G \rightarrow \Z_{2}\times  \Z_{3} \rightarrow 1$ does not split, 
so we have to decompose $C^*(G)$ as a  twisted crossed product of $C^*(H)$ by $\Z_{2}\times \Z_{3}$ and use results of {J.\ }Packer and {I.\ }Raeburn from \cite{PR}.  Actually, when $H$ is a normal subgroup of a group $G$, we give a criterion ensuring  that primitivity of $C^*(H)$ passes over to $C^*(G)$ (see Theorem 1),  and use it to deduce Theorem 2. \footnote{In a forthcoming paper, we will use this criterion to show that $C^*(G)$ is primitive whenever $G$ is the  free product of two nontrivial {\it amenable} groups where at least one of them has more than two elements. This requires an argumentation which is  combinatorially much more involved than in the case of the modular group. }

 \bigskip Murphy mentions in \cite{Mur} that he knows of no example of an icc group  whose full group  C$^*$-algebra is not primitive, but that it is unlikely that such groups do not exist. Now it is almost immediate (cf.\ Proposition 3) that $C^*(G)$ is not primitive whenever $G$ is a nontrivial group having Kazhdan's property (T). As there are many nontrivial icc groups  having property (T), such as $G= PSL(n, \Z)$ for any integer $n \geq 3$ (see \cite{BHV}),  this confirms that the full group  C$^*$-algebra of an icc group is not necessarily primitive. Moreover, as it is known that $PSL(n, \Z)$ is always C$^*$-simple 
(see \cite{BCH, BCH2}), this also illustrates that C$^*$-simplicity of $G$ does not imply that $C^*(G)$ is primitive.

\newpage

\section{On primitivity of full group C$^*$-algebras 
and the modular group}

We use standard notation and terminology in operator algebras; see for example  \cite{Di, Ped, Da}. 
 All Hilbert spaces are assumed to be complex. By a representation of a C$^*$-algebra $A$, we always mean a $*$-homomorphism $\pi : A \to \B(\H)$ into the bounded operators $\B(\H)$ on some Hilbert space $\H$. 
We use the same symbol $\simeq$ to denote unitary equivalence of operators on Hilbert spaces, (unitary) equivalence of representations of a C$^*$-algebra and $*$-isomorphism between C$^*$-algebras.  

\medskip All the groups we consider are assumed to be countable and discrete. If $G$ is such a group, we let  $e_G$, or just $e$, denote its unit. 
 If $G$ acts on a nonempty set $X$ and  $x \in X$, then we say that $x$ is a {\it free point} (for the action of $G$)  whenever
$g \cdot x \neq x$ for all $g \in G, g \neq e$.

 \bigskip 

 Let  $A$ be a  separable C$^*$-algebra and $\widehat{A}$  denote the set of (unitary) equivalence classes of nonzero irreducible representations of $A$. Set
$$\widehat{A}^{\, o}= \{ \, [\pi] \in \widehat{A} \ | \ \pi \ \text{is faithful} \, \} \, . $$
This set is clearly well-defined, and nonempty if and only if $A$ is primitive.

\bigskip 
Assume now that a group $G$ has a normal subgroup $H$ such that $A=C^*(H)$ is primitive and set $K=G/H$. 
 Then   $K$ acts on  $\widehat{A}^{\, o}$ in a natural way.

   \medskip To see this, let $n:K \to G$ be a normalized section for the canonical 
homomorphism $p$ from $G$ onto $K$ (so $n(e_K) = e_G$ and $p \circ n $ gives the identity map on $K$). 

   \medskip Let $\alpha : K \to \text{Aut}(A)$ and $u : K \times K \to A$
 be determined by
$$\alpha_{k}(i_{H}(h))= i_{H}(\, n(k)\, h \,n(k)^{-1}\, ), \ \ \ k \in K, \, h \in H,$$
$$u(k,l)=i_{H}(\, n(k)\, n(l)\, n(kl)^{-1}), \ \ \ k, l \in K,$$
where $i_{H}$ denotes the canonical injection of $H$ into $A$.

   Then $(\alpha, u)$ is a twisted action of $K$ on $A$
  (see \cite{PR} or the Appendix);
 especially, we have
$$\alpha_{k}\, \alpha_{l}= \text{Ad}(u(k, l)) \, \alpha_{kl}, \ k, l \in K,$$
where, as usual,  $\text{Ad}(v)$ denotes the inner automorphism implemented by some unitary $v$ in $A$. 

\medskip This twisted action $(\alpha, u)$ clearly induces an action of $K$ on $\widehat{A}$ given by
$$k \cdot [\pi]= [ \pi \circ \alpha_{k^{-1}}].$$ 
By restriction, we get the {\it natural action} of $K$ on  $\widehat{A}^{o}$, which is easily seen to be independent of the choice of normalized section $n$ for $p$.

\medskip The following result holds:

   \bigskip  
{\bf Theorem 1.} {\it \quad Assume that a group $G$ has a normal subgroup $H$ such that
\vspace{-1ex} \begin{itemize}
\item[a)] $A=C^*(H)$ is primitive, \vspace{-1ex}
\item[b)] $K=G/H$ is amenable,  \vspace{-1ex}
\item[c)] the natural action of $K$ on $\widehat{A}^{o}$  has a free point. 
\end{itemize}
Then $C^*(G)$ is primitive. }

\bigskip

 {\it Proof}. \  We use the notation introduced above and recall that Packer and Raeburn have shown (see \cite[Theorem 4.1]{PR}) that $C^*(G)$ may be decomposed as the  twisted crossed product  associated with $(\alpha, u)$: 
$$C^*(G) \simeq A \times_{\alpha, \, u} K\, .$$
Let   $[\pi] \in \widehat{A}^{o}$ be a free point for the natural  action of $K$. This means that 
$$\pi \circ \alpha_{k}\not \simeq \pi \ \text{for all} \ k \in K, k \neq e\, . $$
Now, this condition implies that
the  induced regular representation Ind $\pi$ of  $A \times_{\alpha, \, u} K$ 
 is irreducible. 
Indeed, as $G$ is discrete, this could be deduced from \cite{Ma}  (see the discussion in \cite[Introduction]{Q1}; see also \cite{Ma2, Ma3, Q2}). For completeness, we give a proof in the Appendix (cf.\  Corollary 6 a)).

 \medskip Further, as $K$ is amenable,  \cite[Theorem 3.1]{PR} gives that Ind $\pi$ is faithful. 
   Altogether, it follows that $C^*(G)$ has a faithful, irreducible
representation, as desired. 

\smallskip 
\hfill $\square$

\bigskip{\bf Remark}. Assume that $G$ has a normal subgroup $H$ and $K=G/H$. 
It would be interesting to find more general conditions than those given in Theorem 1 ensuring that $C^*(G)$ is primitive. However, even for the case where $G$ is the direct product of $H$ and $K$, this is a nontrivial problem. Murphy has shown in \cite[Theorem 2.5]{Mur} that $C^*(H \times K)$ is primitive whenever $C^*(H)$ is primitive and $K$ is amenable and icc. But when for example $\F$ is a free nonabelian group, it is unknown whether  
$C^*(\F \times \F)$ is  primitive or not.
Note that if it should happen that $C^*(\F \times \F)$ is  {\it not} primitive, this would imply that 
$$C^*(\F)\otimes_{max}C^*(\F)\, \,  \not \simeq \, \, C^*(\F)\otimes_{min}C^*(\F).$$ 
Thus, when $\F$ has infinitely many generators, this would solve negatively an open problem of 
E. Kirchberg, which is known to be equivalent to Connes'  famous embedding problem (see \cite{Ki}). 

\endproof

\medskip \bigskip
{\bf Theorem 2.}  {\it \quad  Set $G=PSL(2,\Z)$. Then $C^*(G)$ is primitive. Moreover, there exists an uncountable family of  pairwise inequivalent,  irreducible faithful representations of $C^*(G)$.

} 

\bigskip {\it Proof.} \ Write  $G= \Z_{2}*\Z_{3}= \langle a, b \ | \ a^2 = b^3 = 1 \rangle$ and let $H$ denote the kernel of the canonical homomorphism $p$ from $G$ onto 
$K=\Z_{2}\times \Z_{3} \ ( \, \simeq \Z_{6}\, )$.

\medskip Then   
$H$ is  freely generated as a group  by  
$x_{1}=abab^2$ and $\, x_{2}=ab^2ab$ \\
(see e.g.\  \cite[I.1.3, Proposition 4]{Se}).

   \medskip Set $A=C^*(H)$. Using \cite{Y}, we may pick $[\pi ] \in \widehat{A}^{o}$.    Set 
$$U_{1}=i_{H}(x_{1}),  \, V_{1}= \pi(U_{1}), \ \ U_{2}=i_{H}(x_{2}), \,  V_{2}= \pi(U_{2}), $$
so $V_{1}, V_{2}$ are unitary operators on the separable Hilbert space $\H_{\pi}$ on which $\pi$ acts.    
As shown in the proof \cite[Theorem 6]{Cho}, we may and do assume that  
$V_{2}$ is diagonal relative to some orthonormal basis for  $\H_{\pi}$, 
with (distinct) diagonal entries given by  some $ \mu_{j} \in \T\, ,\,  j \in \N .$

\medskip

For each $\lambda \in \T$,  let $\gamma_{\lambda}$ be the $*$-automorphism of $A$ determined by  
$$\gamma_{\lambda} (U_1)= U_1, \ \gamma_{\lambda} (U_2)= \lambda \, U_2\, ,$$
and set $\pi_{\lambda}= \pi \circ \gamma_{\lambda}$. Clearly, $[\pi_{\lambda}] \in \widehat{A}^{o}$.

\bigskip   We will show that we can pick $\lambda \in \T$ such that $[\pi_{\lambda}]$ is a free point for the natural action of $K$ on 
$\widehat{A}^{o}$. 
 As $K$ is amenable,
the primitivity of $C^*(G)$ will then follow from Theorem 1. 
To pick $\lambda$, we proceed as follows. 

\medskip

As a normalized section for $p:G \to K$, we  choose $n : K \to G$ given by 
$$n(i, j)= a^{i} \, b^{j}, \ \ i \in \{ 0,1\}, \, j \in \{ 0,1,2\} .$$
For each $k=(i,j)\in K$ we let $\alpha_{k}$ 
be the $*$-automorphism of $A$  
used to define the natural action of $K$ on $\widehat{A}^{o}$.

\medskip It is clear that  $[\pi_{\lambda}]$ will be a free point for this action of $K$ if
for each $k \in K, k\neq (0,0),$  we have
$$(\pi_{\lambda}\circ \alpha_{k}) (U_{r}) \not \simeq \pi_{\lambda} (U_{r}) \
 \text{for} \ r= 1 \ \text{or} \ r=2.$$
Some elementary computations give: 
$$\pi_{\lambda} (U_{1})=V_{1}, \, \pi_{\lambda} (U_{2})= \lambda V_{2};$$
 when $ k=(0,1) :  \quad  (\pi_{\lambda}\circ \alpha_{k})(U_{2})= V_{1}^*  ;$ 

\medskip  when $ k=(0,2) :  \quad  (\pi_{\lambda}\circ \alpha_{k})(U_{1})= (\lambda \, V_{2})^*  ;$ 

\medskip  when $ k=(1,0) :  \quad  (\pi_{\lambda}\circ \alpha_{k})(U_{2})= (\lambda  \, V_{2})^*  ;$ 

\medskip  when $ k=(1,1) :  \quad  (\pi_{\lambda}\circ \alpha_{k})(U_{2})= V_{1}  ;$ 

\medskip  when $ k=(1,2) :  \quad  (\pi_{\lambda}\circ \alpha_{k})(U_{1})= \lambda  \, V_{2} \, .$ 

   \medskip It follows that  $[\pi_{\lambda}]$ will be a free point whenever 
$$(*) \quad  V_{1}\not \simeq   \,  \lambda \,V_{2}, \  V_{1} \not \simeq  (\lambda  \, V_{2})^*, \ \lambda \,V_{2} \not \simeq  (\lambda  \, V_{2})^* .$$

Define $\Omega_{1}= \{ \, \lambda \in \T \ | \ V_{1} \simeq  \lambda  \,  V_{2}\, \} $,
 $\quad  \Omega_{2} = \{ \, \lambda \in \T \ | \ V_{1} \simeq  (\lambda  \, V_{2})^*\, \} $,

\smallskip and $ \quad \Omega_{3}= \{ \, \lambda \in \T \ |  \ \lambda V_{2} \simeq  (\lambda  \, V_{2})^*\, \} $.

\medskip As the point spectrum  of $V_{2}$ is  given by $\sigma_{p}(V_{2}) = \{ \mu_{j} \  |Ê \, \, j  \in \N \} \subseteq \T$,  the sets $ \Omega_{1}, \,  \Omega_{2} $ and $\Omega_{3}$ are all countable. 

\medskip Indeed, if  $\Omega_{1}$ was uncountable, then, 
as $\sigma_{p}(V_{1}) = \lambda  \,  \sigma_{p}(V_{2}) $ for all $\lambda \in \Omega_{1}$, 
$\sigma_{p}(V_{1})$ would also be uncountable; as $\H_{\pi}$ is separable, this is impossible.    

In the same way, we see
that $\Omega_{2} $ must be countable.    Finally, if $\Omega _{3}$ was uncountable, then  the equality  
 $$ \ \lambda  \,  \{ \mu_{j} \ |Ê\, \,Ê j \in \N \} = 
 \bar{\lambda} \,  \{ \bar{\mu_{j}} \ |Ê\, \,Ê j \in \N \}$$ would hold for uncountably many $\lambda$'s in $\T$,
 and this is easily seen to be impossible. 

   \medskip Hence, the set $\Omega= \Omega_{1} \cup \Omega_{2} \cup \Omega_{3}$ is countable. Especially, $\Omega \neq \T$ and  $(*)$ holds
for every $\lambda$ in the complement $\Omega^c$ of $\Omega$ in $\T$. Thus, we have shown that $C^*(G)$ is primitive.

\medskip
To prove the second assertion, we consider $\lambda, \lambda' \in \Omega^c$, so $\text{Ind} \, \pi_{\lambda}$ and $\text{Ind} \, \pi_{\lambda'}$ are irreducible and faithful.
A well known argument (adapted to our twisted setting; see Corollary 6 b) in the Appendix) gives that 
$\text{Ind} \, \pi_{\lambda}$ and $\text{Ind} \, \pi_{\lambda'}$ will be inequivalent whenever 
$$\pi_{\lambda} \circ \alpha_j \not \simeq \pi_{\lambda'} \quad \text{ for all} \, \, j \in K\, .$$
Using our previous computations, we see that this will hold whenever 
$$V_{1}\not \simeq   \,  \lambda \,V_{2}\, , \  V_{1} \not \simeq  (\lambda  \, V_{2})^*\, ,$$
$$V_{1}\not \simeq   \,  \lambda' \, V_{2}\, , \  V_{1} \not \simeq  (\lambda'  \, V_{2})^* \, ,$$
$$  \lambda \, V_{2} \not \simeq  \lambda'  \, V_{2}\, , \,  (\lambda V_{2})^* \not \simeq  \lambda'  \, V_{2}\, .$$
The first four conditions are satisfied since $\lambda, \lambda' \in \Omega^c$. 
Set
$$\Omega_{\lambda}= \{Ê\omega \in \T \, | \,  \lambda \, V_{2} \simeq  \omega  \, V_{2}\, \, \text{or} \, \,  (\lambda V_{2})^* \simeq \omega  \, V_{2} \}\, .$$
Then $\Omega_{\lambda}$ is countable (arguing as in the first part of the proof), so $\Omega \cup \Omega_{\lambda}$ is countable.  Hence, if we assume, as we may, that  $\lambda' \in (\Omega \cup \Omega_{\lambda})^c$, then all six conditions above are satisfied, and it follows that 
$\text{Ind} \, \pi_{\lambda}$ and $\text{Ind} \, \pi_{\lambda'}$ are  inequivalent, irreducible and faithful. 

\medskip Proceeding inductively, we may produce in this way a countably infinite family of  pairwise inequivalent,  irreducible faithful representations of $C^*(G)$. In fact, even an uncountable family of such representations does exist. Indeed, observe that  $\text{Ind} \, \pi_{\lambda}$ is an {\it essential} representation of $C^*(G)$, that is, its range contains no compact operators other than zero: 
otherwise,  the  irreducible representations $\text{Ind} \, \pi_{\lambda}$ and $\text{Ind} \, \pi_{\lambda'}$ would have to be equivalent since they have the same kernel (cf.\  \cite[Cor.Ê 4.1.10]{Di}). As $C^*(G)$ is separable, the claim then follows from \cite[Compl\' ementsÊ 4.7.2]{Di}.  

\endproof 

\bigskip
{\bf Remark}. 
Let $G=PSL(2, \Z)$. 
As we have seen in the above proof, $C^*(G)$ has a faithful irreducible representation which is essential. Hence, $C^*(G)$ is antiliminary (cf.\ \cite[Compl\' ements 9.5.4]{Di}). Since $C^*(G)$ is also primitive (and therefore prime), it follows that the pure state space of $C^*(G)$ is weak$^*$ dense in the state space of $C^*(G)$ (cf.\ \cite[Lemme 11.2.4]{Di}). This is also true when $G$ is a nonabelian free group; in fact, this is precisely what Yoshizawa proves in \cite{Y} when $G=\F_2$.

\bigskip Our next observation is quite obvious and surely known to specialists. 

\bigskip 
{\bf Proposition 3.} {\it \quad Let $G$ be a group with Kazhdan's property (T) (see e.g.\ \cite{BHV}) and assume that $C^*(G)$ is primitive. Then $G$ is trivial.}

\bigskip
{\it Proof}. \ Set $A=C^*(G)$. We endow the primitive ideal space Prim$(A)$ of $A$ with its Jacobson (hull-kernel) topology and  $\widehat{A}$ with the weakest topology making the canonical map from $\hat{A}$ onto Prim$(A)$ continuous.   
Since $A$ is primitive, we may pick  $[\pi_0] \in {\widehat{A}}^{\,  o}$. 
As $\{0\}$ is dense in Prim$(A)$, $\{ [\pi_0]\}$ is dense in $\widehat{A}\, $. 

\medskip Let now $\pi_1$ denote the representation of $A$ associated with the trivial one-dimensional unitary representation of $G$. Property (T) means that $[\pi_1]$ is isolated in $\widehat{A}$, i.e.\ 
$\{ [\pi_1]\}$ is  open in $\widehat{A}$. Thus we must have $[\pi_1]=[\pi_0]$. Especially,  $\pi_1$ must be faithful, which implies that $G$ is trivial.
    
\endproof

\bigskip
{\bf Corollary 4.} {\it \quad Set $G=PSL(n, \Z), \, n \geq 3$. Then $G$ is icc, but $C^*(G)$ is not primitive}. 

\bigskip 
{\it Proof}. \ As it is well known that $G$ is icc and has property (T) (see \cite{BHV}), this follows from Proposition 3. 

\endproof

\medskip 
Moreover, as $PSL(n, \Z)$ is always C$^*$-simple 
(cf.\ \cite{BCH, BCH2}), this result also shows that C$^*$-simplicity of  a group $G$ does not imply that $C^*(G)$ is primitive.

\section{Appendix}

We prove here a couple results about induced representations of discrete twisted crossed products, which we could not find explicitely in the literature in the form needed for our purposes.

\bigskip Let  $ (A, K , \alpha,u)$ be a  twisted 
$C^*$-dynamical system as considered by Packer and Raeburn \cite{PR},  where
$A$ is a unital $C^*$-algebra,
$K$ is a discrete group 
with unit $e$
and $(\alpha,u)$ is a twisted action of $K$ on $A$; this means that
 $\alpha$ is a map from $K$ into ${\rm Aut}(A)$, 
the group of $*$-automorphisms of $A$, 
and  $u$ is a map from $K \times K$ into ${\mathcal U}(A)$,
the unitary group of $A$, satisfying
\begin{align*}
\alpha_k \, \alpha_l & = {\rm Ad}(u(k,l)) \, \alpha_{kl} \\
u(k,l) \, u(kl,m) & = \alpha_k(u(l,m)) \, u(k,lm) \\
u(k,e) & = u(e,k) = 1 \ , 
\end{align*}
for all $k, l, m \in K$. (To avoid technicalities, we assume that $A$ is unital; otherwise, one has to assume that the 2-cocycle $u$ takes value in the multiplier algebra of $A$).

\medskip 
The full twisted crossed product $A \times_{\alpha,u} K$ may then be considered as the enveloping
 C$^*$-algebra
of the Banach $*$-algebra $\ell^1(A, K , \alpha,u)$, which consists of the Banach space $ \ell^1(K,A)$  equipped  with product and involution given by
$$( f *g ) (l) = \sum_{k\in K} f(k) \alpha_k(g(k^{-1}l)) u(k, k^{-1}l)$$
$$f^{*}(l) = u(l, l^{-1})^{*} \alpha_l (f(l^{-1}))^{*} $$
$f, g \in \ell^1(K,A), \, l \in K$.

\medskip We let $i_K$ and $i_A$ denote the canonical injections of $K$ and $A$ into $A \times_{\alpha,u} K$, respectively.

\medskip
Let now $\pi$ be a nondegenerate representation of $A$
on some Hilbert space $\H=\H_\pi$ and let  $\pi_\alpha$ be the associated representation of $A$ on 
$\H_K=\ell^2(K,\H)$ 
defined by
$$(\pi_\alpha(a) \xi)(k) = \pi(\alpha_{k^{-1}}(a))\, \xi(k) \ , 
\, a \in A, \xi \in \H_K, k \in K \, .$$ 
For every $k \in K$, let $\lambda_u(k)$ be the unitary operator on 
$\H_K$ given by
$$(\lambda_u(k)\xi)(l) = \pi(u(l^{-1},k))\, \xi(k^{-1}l), \, k, l \in K, \,  \xi \in 
\H_K\,. $$
The pair $(\pi_\alpha,\lambda_u)$ is then a covariant representation 
of $(A,K,\alpha,u)$, that is,
\begin{align*}
\pi_\alpha(\alpha_k(a)) & = {\rm Ad}(\lambda_u(k))(\pi_\alpha(a))\\
\lambda_u(k) \lambda_u(l) & = \pi_\alpha(u(k,l)) \lambda_u(kl) 
\end{align*}
for all $k,l \in K$ and $a \in A$. (Note that we follow  \cite{ZM} here, while the "right" version is used in \cite{PR}).

\medskip This covariant representation  induces  a nondegenerate representation ${\rm Ind}\,  \pi$ of  $A \times_{\alpha,u} K$ on $\H_K$ determined by 
$$({\rm Ind} \, \pi)(f)= \sum_{k \in K}\pi_\alpha(f(k)) \lambda_u(k), \quad  f \in \ell^1(K, A)\, ,$$
that is, by 
$$({\rm Ind} \, \pi)(i_A(a)) = \pi_\alpha(a) \, , \, ({\rm Ind} \, \pi)(i_K(k)) = \lambda_u(k)\, \, , \, a \in A\, , k \in K\, .$$

\newpage For each $k \in K$,  let $\H_k$ denote the  copy of $\H$ in $\H_K$ given by   $$\H_k=\{ \xi \in \H_K\, | \, \xi(l) = 0 \, \, \text{for all} \, \, l \in K\, , \, l\neq k\}\, ,$$ giving us the  natural direct sum  decomposition $\H_K = \oplus_{k \in K} \H_{k}$. 

\medskip Assume now that $\pi'$  is a nondegenerate representation of $A$ on $\H'$ and denote by $(\pi'_\alpha,\lambda'_u)$ the associated covariant representation 
of $(A,K,\alpha,u)$ on $\H'_{K}$. 

\medskip Let $T \in \B(\H_K, \H'_{K})$. Denote by $[T_{k,l}]_{k, l \in K}$ the matrix of $T$ with respect to the natural direct sum  decompositions of $\H_K$ and $\H'_K$, and identify each $T_{k,l}$ as an element in $\B(\H, \H')$.  

\medskip Hence, if $\eta \in \H$ and $ k, l \in K$, then $T_{k,l}\,  \eta = (T\,  \eta_l)(k)$, where $\eta_l \in \H_K$ is given by 
  $\eta_l(k) = \eta$ when $k=l$, and $\eta_l(k) = 0$ otherwise.
 
\medskip Some tedious (but straightforward) computations give:

\medskip $(1)\quad   (T\, \pi_\alpha(a))_{k,l}= T_{k,l} \,  \pi(\alpha_{l^{-1}}(a)) \, , \quad  \, \, \, (\pi'_\alpha(a)\, T)_{k,l}= \pi'(\alpha_{k^{-1}}(a))\, T_{k,l} \, \, ,$

\medskip $(2) \ (T\, \lambda_u(j))_{k,l}= T_{k, jl}\, \pi(u(l^{-1}j^{-1}, j))  \, , \,  (\lambda'_u(j)\, T)_{k,l}= \pi'(u(k^{-1}, j))\, T_{j^{-1}k, l} \,.$

 \bigskip {\bf Proposition 5.} {\it 
 
 \medskip Assume $\pi$ and $\pi'$ are irreducible, and $\pi \circ \alpha_j \not \simeq \pi'$ for all
$j \in K, \, j \neq e$. 

\smallskip 
Let $T \in \B(\H_K, \H'_K)$  intertwine ${\rm Ind} \, \pi$ and ${\rm Ind} \, \pi'$.

\medskip Then $T_{k,k}$ intertwines $\pi$ and $\pi'$ for all $k \in K$. Further, $T$ is decomposable, that is, $T_{k,l} = 0$ for all $k \neq l$ in $K$}.

\bigskip
{\it Proof}. 
We first note that $ \ T\, \pi_\alpha(a) =  \pi'_\alpha(a)\, T$ for all $a \in A$.
Using (1), we then get 

\medskip $(3) \quad T_{k,l}\,  \pi(\alpha_{l^{-1}}(a)) =   \pi' (\alpha_{k^{-1}}(a)) \,  T_{k,l} \, \, $ for all $\,  k,l \in K, \, a \in A.$

\medskip Letting $l=k$, this clearly implies that $T_{k,k}$ intertwines $\pi$ and $\pi'$ for all $k \in K$. 

\medskip 
Assume now that $k\neq l$. Using (3) with $a=\alpha_k(b)$, we get 

\medskip 
$(4) \quad  T_{k,l}\,  (\pi \circ {\rm Ad}(u(l^{-1}, k))  \circ \alpha_{l^{-1}k})(b) = (\pi' \circ  {\rm Ad}(u(k^{-1}, k)))(b) \, T_{k,l}$

\medskip for all $b \in A$. 

\medskip From the assumption, we have  $\pi' \not \simeq \pi \circ \alpha_{l^{-1}k}$. Hence, it follows that $\pi \circ  {\rm Ad}(u(l^{-1}, k) \circ \alpha_{l^{-1}k}$
 and $\pi' \circ {\rm Ad}(u(k^{-1}, k))$ are irreducible and inequivalent. But (4) says that $T_{k,l}$ intertwines these two representations of $A$, and we can therefore conclude that $T_{k,l} = 0$.

\endproof

\bigskip The following corollary is due
 to Zeller-Meier  in the case where 
 $u$  takes values in the center  of $A$ (see  \cite[Propositions 3.8 and 4.4]{ZM}). Part a) could be deduced from \cite[Theorem]{Q2},
but as we also need part b), we prove both. 

\bigskip

{\bf Corollary 6.} {\it 

 \medskip a) $ {\rm Ind} \, \pi$ is irreducible whenever $\pi$ is irreducible and the stabilizer subgroup $K_\pi = \{k\in K \, | \, \pi \circ \alpha_k \simeq \pi \}$ is trivial.

\medskip b) Assume that $\pi$ and $\pi'$ both are irreducible.

 \medskip Then   $\, {\rm Ind}  \, \pi \, \not \simeq \, {\rm Ind} \, \pi'$ whenever $\pi \circ \alpha_j  \, \not \simeq \, \pi'$ for all $j \in K$.
}

\bigskip 
{\it Proof}.  
a) Suppose that  $\pi$ is irreducible and  $K_\pi $ is trivial.

\medskip Let  $T \in \B(\H_K)$ lie in the commutant of $({\rm Ind} \, \pi)( A\times_{\alpha,u} K)$.

\medskip Using Proposition 5  with $\pi'=\pi$, it follows that $T$ is decomposable and
 $T_{k,k} \in \pi(A)'$ for all $k \in K$. As  $\pi$ is irreducible, this gives that
$T_{k,k} \in \C \, I_\H$ for all $k \in K$.

\medskip Further, we have $  T \, \lambda_u(j) =  \lambda_u(j)\, T$ for all $j \in K$.

 \medskip Using this and  (2),  we get 
$$\pi(u(k^{-1}, k l^{-1}))\, T_{k,k} = T_{k,k}\, \pi(u(k^{-1}, k l^{-1})) = (T\, \lambda_u(k l^{-1}))_{k,l}$$ 
$$= (\lambda_u(k l^{-1}) \, T)_{k,l} = \pi(u(k^{-1}, k l^{-1}))\, T_{l,l}\, , $$
which  implies that $T_{k,k}=T_{l,l}$  for all $k,l \in K$.

\medskip Altogether, this means that $T$ is a scalar multiple of the identity operator on $\H_K$. Hence we have shown that ${\rm Ind} \, \pi$ is irreducible, as desired. 

\bigskip b) Assume that $\pi$ and $\pi'$ both are irreducible and $\pi \circ \alpha_j  \, \not \simeq \, \pi'$ for all $j \in K$.

\medskip Let  $T \in \B(\H_K, \H'_K)$  intertwine ${\rm Ind} \, \pi$ and ${\rm Ind} \, \pi'$. It follows from Proposition 5 that
$T_{k,l} = 0$ for all $k,l \in K, \, k\neq l$, and that $T_{k,k}$ intertwine $\pi$ and $\pi'$ for all $k \in K$. As $\pi \not \simeq \pi'$ by assumption,  we also have $T_{k,k}=0$ for all $k\in K$. Hence, $T=0$. This shows that  ${\rm Ind}  \, \pi \, \not \simeq \, {\rm Ind} \, \pi'$, as desired.

 \endproof
 
Actually, both implications converse to those stated in a) and b) of Corollary 6 also hold (as in \cite{ZM}). However, since we don't need these in this paper, we skip the proofs.

\end{flushleft}

 \end{document}